\documentclass[12pt]{amsart}
\usepackage{amsmath}
\usepackage{amssymb}
\usepackage{amsthm}
\newcommand{\ignore}[1]{\relax}

\newcommand{\C}{\mathbb C}
\newcommand{\R}{\mathbb R}
\newcommand{\Z}{\mathbb Z}

\newcommand{\U}{\mathcal U}

\newcommand{\PP}{\mathbb P}

\theoremstyle{definition}

\theoremstyle{remark}

\newcommand{\ev}{\operatorname{ev}}

\newcommand{\T}{\mathbb{T}}

\renewcommand{\setminus}{\smallsetminus}

\newcommand{\MM}{{\mathcal M}}
\newcommand{\Jac}{\operatorname{Jac}}
\newcommand{\Div}{\operatorname{Div}}
\newcommand{\Pic}{\operatorname{Pic}}

\newcommand{\OO}{\mathcal O}

\begin{document}

\title
[Introduction to Tropical Geometry]
{Introduction to Tropical Geometry,\\
notes from the IMPA lectures, Summer 2007}
\author{Grigory Mikhalkin}
\thanks{The author is partially supported by the Canada Research
Chair program, NSERC and PREA}
\address{Department of Mathematics\\
University of Toronto\\
Toronto, Ontario M5S 2E4 Canada} \email{mikha@math.toronto.edu}
\maketitle

\section{Tropical varieties}
\subsection{Generalities}
\subsubsection*{\bf A}
In our viewpoint an algebraic variety is a pair $(X,\OO)$ consisting
of a topological space $X$ and a sheaf of function (to the base field) called
{\em regular functions} subject to some nice geometric properties.

\subsubsection*{\bf B} In tropical geometry the base field is not strictly speaking a field,
but only a semifield $\T$. It consists of the real numbers augmented by $-\infty$
and enhanced by two arithmetic operation: $``x+y"=\max\{x,y\}$ and
$``xy"=x+y$. These operations are commutative associate and satisfy to
the distribution law. We have the additive zero $0_{\T}=-\infty$
and the multiplicative unit $1_{\T}=0$. The multiplicative group $\T^{\times}=
\T\setminus\{-\infty\}$
is an honest Abelian group with a division, but there is no subtraction
whatsoever. The tropical addition is idempotent, $``x+x"=x$.

\subsubsection*{\bf C}
A tropical Laurent polynomial in $n$ variables is a function
$$f(x_1,\dots,x_n)=``\sum a_{j_1\dots j_n}x_1^{j_1}\dots x_n^{j_n}\ "
=\max\limits_j j_1x_1+\dots+j_nx_n,$$
where the sum is finite and the indices $j_l$ are integers.
The space $(\T^{\times})^approx\R^n$ equipped with the Euclidean
topology and considered with the sheaf of tropical Laurent polynomial
is an example of a tropical variety (analogous to the torus $(\C^*)^n$).

\subsubsection*{\bf D}
One can evaluate a Laurent polynomial at any point of $\R^n$
and at some points of $\T^n\supset\R^n$. Namely, a Laurent polynomial
is regular at a point of $\T^n$ if it has a limit at this point bounded from
above (possibly $-\infty$ but not $+\infty$). This defines a regular
sheaf on the {\em tropical affine space} $\T^n$. The only functions regular
on the entire $\T^n$ are honest polynomials (such that all the powers $j_l$
are non-negative).

\subsubsection*{\bf E} Starting from $\R^n$ we can get compact tropical varieties
by taking a quotient under an arbitrary discrete lattice $\Lambda\approx\Z^n$
embedded to $\R^n$. Later we'll see that some of these tropical tori are
algebraic while some are not (as it is the case for abelian varieties and
other complex tori). This gives us some examples of tropical varieties modeled
on manifolds, but the underlying topological space for most tropical varieties
is not a manifold but only a polyhedral complex.

\subsection{Polyhedral complexes and tropical varieties}
\subsubsection*{\bf A} A convex polyhedron $P$ with integer slope in $\R^N$
is an intersection of half-spaces in $\R^N$ defined by equations
$ax\ge b$. Here $x\in\R^N$, $a\in\Z^N$, $b\in\R$.
A convex polyhedron in $\T^N$ is a closure of $P$ in $\T^N\supset\R^N$.
A {\em dimension} of such a polyhedron is its usual topological dimension.

\subsubsection*{\bf B} An (integer-slope $n$-dimensional) polyhedral complex $P$
in $\R^N$ is a union
of a collection of convex polyhedra with integer slope of dimension $n$ such that the
intersection of any subcollection is a common face. Again, to get a polyhedral
complex in the tropical affine space $\T^N$ we simply take the closure.

\subsubsection*{\bf C} Sedentarity of a point in $\T^N$ is the number of
coordinates equal to $-\infty$. The sedentarity of a face of $P$
is the sedentarity of a generic point there. We call points
and faces of zero sedentarity {\em inner}.

\subsubsection*{\bf D}\label{bal}
Consider an $(n-1)$-dimensional inner face $E$ of $P$. Let $F_1,
\dots, F_k$ be the facets ($n$-dimensional faces) adjacent to $E$.
Lel $\lambda:\R^N\to\R^{N-n+1}$ be the projection whose kernel is parallel to $E$.
The facets $F_1,\dots,F_k$ project to rays adjacent to $\lambda(E)$.
Clearly, these rays have integer slope in $\R^{N-n+1}$.
Let $v_1,\dots,v_k\in\Z^{N-n+1}$ be the primitive integer vectors parallel
to $\lambda(F_1),\dots,\lambda(F_k)$. We say that $P$ is {\em balanced}
if $\sum\limits_{j=1}^kv_j=0$.

\subsubsection*{\bf E} An $n$-dimensional tropical variety is a pair $(X,\OO)$
if $X$ is a Hausdorff space and can be covered by finitely many open sets
$V_\alpha\subset X$ such that for each $\alpha$ there exists an open set
$U_\alpha\supset V_\alpha$ in $X$, a balanced polyhedral
complex $P_\alpha\subset\T^{N_\alpha}$ and an embedding
$\iota_\alpha:\U_\alpha\to P_\alpha$
with the property that the closure of $\iota_\alpha(V_\alpha)$ is contained
in $\iota_\alpha(U_\alpha)$ and that the restriction of the sheaf $\OO$ to $U_\alpha$
coincides with the restriction of $\OO_{\T^{N_\alpha}}$ under $\iota_\alpha$.

\subsubsection*{\bf F}
It is easy to see that the product of a finite number of
balanced polyhedral complexes is again a balanced polyhedral complex and
thus the product of a finite number of
tropical varieties is again a tropical variety.

\subsubsection*{\bf G} {\em Remark.} It is also possible to define a {\em smooth} tropical
variety by modeling it on a {\em simply balanced polyhedral complex}.
A slightly different class of tropical varieties is formed by images
of smooth tropical varieties under tropical morphisms. This leads
to a more proper definition of tropical varieties. However, for the
sake of simplicity in these lectures we use the definition given above.

\subsection{Hypersurfaces and complete intersections}

\subsubsection*{\bf A} Let $f:\T^N\to \T$ be a regular function (tropical
polynomial). The {\em hypersurface} $V_f$ defined by $f$ is the set of points
in $\T^n$ where $``\frac{1_{\T}}{f}"=-f$ is not regular. In other words,
$V_f$ is the {\em corner locus} of the convex piecewise-linear function $f$.

\subsubsection*{\bf B} The hypersurface $V_f$ is an $(n-1)$-dimensional
polyhedral complex in $\T^N$.
However, in general it is NOT balanced. Nevertheless, one may associate
a natural number (called the {\em weight}) that makes it into a {\em weighted
balanced} polyhedral complex (so that we get zero when we modify the equality
from \ref{bal}.D by weighting all the primitive vectors.

\subsubsection*{\bf C} Similarly, one can define a hypersurface associated
to a regular function $f$ on a tropical variety $X$ (or an open set $U\subset X$).

\subsubsection*{\bf D} The weights on the facets play the role of the orders
of zeroes of $f$. Alternatively they may be defined recursively by saying
that that the order of $``fg"$ at a facet $F$ is $a+b$ if $a$ is the weight
of $f$ at $F$ and $b$ is the weight of $g$ at $F$.
If all the weights of the facets of $V_f$ are 1 then $V_f$ is a tropical subvariety.
Otherwise it is only a tropical {\em (Cartier) divisor}.

\subsubsection*{\bf E} More generally, one can define a tropical {\em $k$-cycle}
$A\subset X$ by defining them as (locally) weighted balanced $k$-dimensional
polyhedral subcomplex inside
the tropical variety. It is convenient to allow not only natural, but integer weights.

\subsubsection*{\bf F}
{\em Warning.} It may happen that $A$ is an $(n-1)$-cycle in a tropical $n$-dimensional
variety $X$, but not a hypersurface for any regular function $f$ (even locally and
even in the case when $X$ is smooth). In other words in tropical geometry there
exist non-Cartier divisors even for smooth varieties.

\subsubsection*{\bf G}
A tropical subvariety $Y\subset X$ is a balanced {\em unweighted}
subcomplex (i.e. such a divisor that all the weights are 1).

\subsubsection*{\bf H}
{\em Example.} A tropical cubic polynomial in 3 variables defines a
smooth hypersurface in $\T^3$ if the (partially defined) coefficient
function is convex. We may find generic cubics with 27 distinct
lines as 1-cycles. ({\em Warning:} for some other generic choices
we have positive-dimensional families of lines, so there have infinitely
many lines, unlike the complex case.)

\subsection{Projective spaces and cycles there}
\subsubsection*{\bf A}
Let $\phi:\T^\times\to\T^\times$ be defined by $z\mapsto ``\frac{1}{z}"$.
The result of gluing of two copies of $\T\supset\T^\times$ along $\phi$
is called the tropical projective line and is denoted with $\T\PP^1$.
This is the tropical counterpart of the Riemann sphere.

\subsubsection*{\bf B}
As a topological space it is just a closed interval.
Geometrically the tropical structure on $\T\PP^1$ may be
reformulated as the {\em integer affine structure} on $\R=\T^\times\subset\T\PP^1$,
the stratum of sedentarity 0.
Indeed, to be able to define tropical polynomials we only need to know
what is an affine-linear function with an integer slope.
The boundary divisor $\T\PP^1\setminus\R^n$ consists of two components,
each isomorphic to the single-point space $\T\PP^0$.

\subsubsection*{\bf C}
In a similar way we may glue $n+1$ copies of $\T^n$ to get the
tropical projective space $\T\PP^n$. The sedentarity 0 stratum is
$\R^n$. Again its (tautological) integer affine-linear structure encodes
the tropical structure on $\T\PP^n$. The complement $\T\PP^n\setminus\R^n$
consists of $(n+1)$ divisors called the boundary divisors
and isomorphic to $\T\PP^{n-1}$.

\subsubsection*{\bf D}
A point has
sedentarity $k$ if and only if it is contained in exactly $k$ different
boundary divisors. Points of sedentarity $k$ are thus divided into
$\begin{pmatrix} n+1\\ k\end{pmatrix}$ strata that correspond to the
$k$-dimensional faces of the (topological) $n$-simplex $\T\PP^n$.
Each individual stratum is isomorphic to $\T\PP^k$.

\subsubsection*{\bf E}
We may define the projective {\em degree} of a $k$-cycle $A\subset\T\PP^n$
by computing its {\em intersection number} with one of the boundary $\T\PP^{n-k}$.
If we equip each intersection point with the proper multiplicity then
(thanks to the balancing condition) it does not depend on the choice of $\T\PP^{n-k}$.

\subsubsection*{\bf F}
Instead of giving the definition of the intersection number in this boundary case (which is
also possible, but more complicated) we may choose
yet another (non-boundary and easier to use) representative of $\T\PP^{n-k}$.
Take the negative parts of the $n$ coordinate axes in $\R^n$
and the ray $(t,\dots,t), t\ge 0$. Alltogether these are $n+1$ rays emanating from
the origin.
For each $(n-k)$ out of these $n+1$
rays we take their convex hull. This gives us
$\begin{pmatrix} n+1\\ n-k\end{pmatrix}$ $(n-k)$-dimensional quadrants (each
isomorphic to $\R^{n-k}_{\ge 0}$.
We denote the union of all the closure of these quadrants in
with $L_{n-k}\subset\T\PP^n$. It is a tropical variety.

\subsubsection*{\bf G}
Note that any translation in $\R^n$ may be extended to $\T\PP^n$ by taking
the closure. We may move $L_{n-k}$ by such translations $\tau$. Note that for
a generic choice of $\tau$ the polyhedral complexes $A$ and $\tau L_{n-k}$
are transverse, i.e. they intersect in a finitely many points of sedentarity 0
and these points are in the interior of the facets of $A$ and $\tau L_{n-k}$
respectively.

\subsubsection*{\bf H}
At a point of a transverse intersection at interior points
of a facet of dimension $k$ and and a facet of dimension $(n-k)$
the definition of the local intersection number is especially easy.
Indeed, locally such an intersection is given by two transverse
affine-linear subspaces. Let $\Lambda_1\subset\Z^n$ (resp. $\Lambda_2$)
be the integer vectors parallel to the first (resp. second) subspace.
We define the local intersection number as the index of the subgroup
$\Lambda_1+\Lambda_2\subset\Z^n$ multiplied by the weights of the
corresponding facets. The balancing condition ensures that the intersection
number is invariant with respect to the deformations, in particular to
translations in $\R^n$ of one of the cycles.

\subsubsection*{\bf I} Note that the subspace $L_{k}$ is not even homeomorphic
to $\T\PP^n$. Nevertheless, it can be considered as a result of deformation
of the boundary $\T\PP^k$. The tropical varieties $L_k$ and $\T\PP^k$
are not isomorphic but they are {\em equivalent}.

\subsection{Equivalence of tropical varieties}

\subsubsection*{\bf A}
Let $f:\T^n\to\T$ be a regular function. It is easy to see that the
set-theoretical graph of $f$ is NOT a hypersurface in $\T^{n+1}=\T^n\times\T$.
Namely, the balancing condition fails at the corner locus of $f$ as there the
function is strictly convex.

\subsubsection*{\bf B} In the same time the set-theoretical graph of $f$
is contained in the hypersurface $\Gamma_f$
defined by $y+f(x)$ (where $x\in\T^n,y\in\T$).
This hypersurface can be considered as the ``complete tropical graph" of $f$.
It is different from the set-theoretical graph by attaching a ray
$\{x\}\times[-\infty,f(x)]$ for each $x\in V_f$. Furthermore,
the weights of the corresponding facets in $V_f$ and $\Gamma_f$ coincide.
We note that $\Gamma_f$ is NOT isomorphic to $\T^n$. We postulate that
$\Gamma_f$ and $\T^n$ are {equivalent}.

\subsubsection*{\bf C} More generally, the same situation holds for
a regular function on a tropical variety $X$ or an open set $U\subset X$.
Recall that the hypersurface $V_f\subset U$ is still well-defined as
the set of points where $``\frac{1_{\T}}{f}"$ is not regular and
we still have weights at the facets of $V_f$. Thus we can still
define $\Gamma_f$ by attaching $\{x\}\times[-\infty,f(x)]$ for each
$x\in V_f$ and enhancing the new facets with the weights coinciding
with the weights of the corresponding facets in $V_f$.
Alternatively one can define the weights of the new facets
of $\Gamma_f$ so as to make $\Gamma_f\subset\T^{n+1}$
to be a weighted balanced polyhedral complex.

\subsubsection*{\bf D} {\em Remark.} Once we'll develop
the tropical intersection theory below we'll have yet
another alternative way of definition of $\Gamma_f$
(together with the weights on its facets). Namely,
$X$ locally coincides with a balanced polyhedral complex
$P\subset\T^N$ and $f$ comes as a restriction of a tropical
Laurent polynomial $F:\T^N\to\T$.
We can define $\Gamma_f$ locally as the {\em stable intersection}
of $P\times\T$ and the hypersurface of $y+F(x)$ in $\T^{N+1}$.

\subsubsection*{\bf E} {\em Definition.}
We say that a map $Y\to X$ between tropical varieties is a
{\em tropical modification} if locally it coincides with
the projection $\Gamma_f\to U$ for some regular $f:U\to\T$.
Two tropical varieties are called {\em equivalent} if they
can be connected with a sequence of tropical modifications
(an their inverse operations).

\subsubsection*{\bf F} {\em Remark.} If $X$ is a tropical variety
then the tropical structure on $Y$ is defined by $f$ if all the
weights of $\Gamma_f$ are equal to 1. We may alternatively define
SMOOTH $n$-dimensional tropical varieties as topological spaces
with sheaves of functions locally equivalent to $\T^n$
after adding some finiteness and separability conditions.
Note that
this definition is analogous to analytic definition of complex manifolds:
{\em a complex manifold is a topological space with a sheaf of functions
locally isomorphic to $\C^n$.}

\subsection{Tropical morphisms}

\subsubsection*{\bf A} A map $\phi:Y\to X$ is called a {\em tropical morphism}
if for every $y\in Y$ there exists a neighborhood $y\in U\subset Y$
and an affine-linear map $\Phi:\R^N\to\R^M$ such that the linear part
of $\Phi$ is defined over $\Z$ and such that under the charts
$\iota:U\subset\T^N$ and $\phi(U)\subset\T^M$ the map $\phi$ is given
by the closure of $\Phi$ to a partially defined map from $\T^N$ to $\T^M$.
Here we demand that whenever $z\in\iota(U)$ the closure of $\Phi$
is actually defined at $z$.

\subsubsection*{\bf B} Note that if $X$ and $\tilde{X}$ are equivalent
tropical varieties then there exists their ``common resolution" $\bar{X}$,
i.e. a tropical variety equivalent to $X$ and $\tilde{X}$ with tropical
morphisms $\phi:\bar{X}\to X$ and $\tilde{\phi}:\bar{X}\to\tilde{X}$
that are both composition of tropical modifications. Indeed, if $\psi:V\to U$
is a tropical modification corresponding to a regular function
$g:U\to\T$ and $f:U\to\T$ is another regular function then
$f\circ\psi$ is tropical on $\psi^{-1}(U)$ and we can make consider
a tropical modification $W\to V$ corresponding to $f\circ\psi$.
Furthermore, we can obtain the same variety
$W$ by making a tropical modification along
$f$ first and then by the pull-back of $g$.
In any case $W$ is obtained from $U$ by attaching rays along the hypersurfaces
of $g$ and $f$ and by attaching a positive quadrant along their
{\em stable intersection} (see below) which is symmetric.

\subsubsection*{\bf C}
Two tropical morphisms $\phi:Y\to X$ and $\tilde\phi:\tilde{Y}\to\tilde{X}$
are called equivalent if $\tilde{Y}$ is equivalent to $Y$, $\tilde{X}$
is equivalent to $X$ and $\tilde{\phi}$ corresponds to $\phi$ under
these equivalences (in other words, they can be lifted to the same morphism
from the common resolution of $Y$ and $\tilde{Y}$ to the common resolution
of $X$ and $\tilde{X}$.

\section{Curves}
\subsection{Tropical curves as metric graphs}

\subsubsection*{\bf A}
Zero-dimensional tropical varieties are just disjoint union of points.
The first interesting example of tropical varieties appears in dimension 1.
Clearly, any tropical curve is a graph. A compact tropical curve is a finite
graph $\bar\Gamma$. Let us assume in addition that $\bar\Gamma$ is connected.
Denote with $\Gamma$ the complement of the set of 1-valent vertices in $\bar\Gamma$.

\subsubsection*{\bf B}
The graph $\bar\Gamma$ has an additional geometric structure that carries the
same information as enhancing it with a structure sheaf of regular tropical functions.
This structure is the complete inner metric on $\Gamma$. Recall that
a {\em leaf} is an edge adjacent to a 1-valent vertex. Completeness of the metric
ensures that all leaves have infinite length. The other edges (called {\em inner edges})
have finite length.

\subsubsection*{\bf C}
Once $\Gamma$ is equipped with a complete inner metric we can reconstruct
the structure sheaf $\OO_{\bar{\Gamma}}$ in the following way. For each point
inside an edge a metric gives a well-defined (up to translation) embedding
to $\T^\times\approx\R$ and we pull back the regular functions (the tropical
Laurent polynomials). Near a 1-valent vertex this embedding can be extended
to an embedding to $\T$. Finally, for a vertex of valence $k+1\ge 3$
we have a standard model: a curve $C_k\subset\R^k$ obtained as a union
of the negative halves of $k$ coordinate axes with the diagonal
$(t,\dots,t)$, $t\ge 0$. The resulting union is balanced, restrictions
of the tropical Laurent polynomials in $k$ variables to $C_k$ locally
form the structure sheaf.

\subsubsection*{\bf D}
Conversely, if $\Gamma$ is a tropical curve then we have the notion of
a {\em primitive tangent vector} to $\Gamma$. At a point inside an edge
such vector is well-defined up to the sign. The metric on $\Gamma$ is defined
by making these vectors unitary. The finiteness condition ensures that every
leaf has an infinite length.

\subsubsection*{\bf E} {\em Remark}
More generally, a tropical structure on a tropical variety can be
rephrased in terms of the so-called {\em integer affine structure}.
While traditionally such structure is considered only on smooth manifolds
we can express the tropical structure as its extension to certain
polyhedral complexes.

\subsection{Tropical modifications for curves}

\subsubsection*{\bf A}
Tropical modifications for curves are especially easy since hypersurfaces
in curves are disjoint unions of points (with some natural weights).
We restrict our attention to the case of weight 1, i.e. a modification
along a simple point. In the language of metric graph such a modification
results in attaching an interval of infinite length to this point.
We easily see that all compact metric trees are equivalent. From now on
we restrict our consideration to the case of compact and connected
tropical curves.

\subsubsection*{\bf B}
If we have a finite collection of marked points of sedentarity 0
then we can make tropical modifications at them to replace the curve
with an equivalent model so that all the marked points have sedentarity 1,
i.e. are the 1-valent vertices (ends) of the curve. Furthermore, if the number
of the marked points is at least 2 then we can do the inverse tropical
modifications so that the ends are precisely the marked points.

\subsection{Differential forms and their integrals}

\subsubsection*{\bf A}
A constant differential $k$-form $\alpha$ on $\R^N$ is just an element
of $\Lambda^k((\R^N)^*)$. We say that it extends to a point
$x\in\T^N\setminus\R^N$ if for any infinite coordinate of $x$
the corresponding basis vector is in the kernel of $\alpha$.

\subsubsection*{\bf B}
A {\em regular} differential $k$-form on a tropical variety $X$
is a differential form that can be locally obtained as a restriction
of a constant differential form in $\T^N$ defined at all points of $X$.

\subsubsection*{\bf C}
We may integrate differential forms against respective chains.
E.g. let
$\alpha$ be a 1-form on a tropical curve $C$ and let $\gamma:[0,1]\to C$
be a path.
Note that because of the extension condition any regular 1-form
on a leaf of a compact tropical curve must vanish.
Thus the integral $\int\limits_\gamma\alpha$ is a real number.

\subsection{The genus of tropical curve and its Jacobian}
\subsubsection*{\bf A}
All regular 1-forms on a compact connected tropical curve $C$
form a finite-dimensional
(real) vector space $\Omega(C)$. Its dimension is equal to $b_1(C)=\dim H_1(C)$.
This number is called the {\em genus} of $C$.

\subsubsection*{\bf B}
Integration along 1-cycles gives a map
$H_1(C;\Z)\to\Omega^*(C)$, this map is an embedding and its
image is called the {\em periods}.
The quotient $$\Jac(C)=\Omega^*(C)/H_1(C;\Z)$$ is called
the {\em Jacobian} of $C$. Topologically it is a torus,
but it is naturally equipped with a tropical structure coming
from the tautological integer affine structure
on $\Omega(C)$ (its integer affine structure is coming from the
regular forms that take integer values on integer tangent vectors to $C$).

\subsection{The divisor group and the Picard group}
\subsubsection*{\bf A}
The Jacobian $\Jac(C)$ can also be interpreted as the Picard group
of degree $0$. Let us recall its definition. The divisor group $\Div$ is defined
as all formal (finite) linear combinations $D=\sum a_jx_j$, where $a_j\in\Z$ and
$x_j\in C$. The {\em degree} of the divisor $D$ is defined as $\sum a_j\in\Z$.
The divisors of degree 0 form a group $\Div_0$.

\subsubsection*{\bf B}
Note that there are no non-constant globally defined regular functions on compact
tropical curves. As in the complex case the reason is the maximum principle:
at the point of its maximum a regular function must be locally constant.
Nevertheless, the function which at every point can be presented as the
locally tropical quotient of two regular functions can be defined globally.
Such functions are called {\em rational functions}, they are piecewise-linear,
though no longer necessarily convex.

\subsubsection*{\bf C}
A point $x\in C$ is called a zero of order $m$ of a rational function $\phi:C\to\T$
if near $x$ the function $\phi$ coincides with a regular function with a zero
of order $m$. Similarly, $x\in C$ is called a pole of order $m$ if near $x$
$``\frac{1}{\phi}"=-\phi$ coincides with a regular function with a zero
of order $m$.

\subsubsection*{\bf D}
The divisor of a rational function $\phi$ is the linear combination
of all its zeroes taken with the coefficients equal to their order and
all its poles taken with the coefficients equal to minus their order.
Such divisors are called {\em principal}. Their degree is always zero.
Two divisors are called {\em linearly equivalent} if their difference is
a principal divisor. The quotient of the divisor group by the linear
equivalence is called the {\em Picard group} and is denoted with $\Pic(C)$.
The quotient of the divisors of degree $d$ by the linear equivalent
is called $\Pic_d(C)$.

\subsection{The Abel-Jacobi theorem}
\subsubsection*{\bf A}
We have a map from the product of $k$ copies of the curve $C$
(which is a tropical variety) to $\Pic_k(C)$. Indeed, any $k$
points form a divisor of degree $k$ and we can take its equivalence
class in $\Pic(C)$. Note that $\Pic_0(C)$ and $\Pic_k(C)$ can be identified
by an invertible map given by adding a copy of a fixed point in $C$ $k$ times.
Thus $\Pic_k(C)\approx\Pic_0(C)$ as a set.

\subsubsection*{\bf B}
Integration gives a natural map from $\Pic_0(C)$ to $\Jac(C)$.
Namely, any element of $\Pic_0(C)$ is a $0$-dimensional cycle in $C$
and thus it bounds a 1-chain $\gamma$. Integration along $\gamma$
gives a linear $\R$-valued functional on the regular 1-forms $\Omega(C)$,
i.e. an element of $\Omega^*(C)$. The 1-chain $\gamma$ is defined up
to 1-cycle in $C$. Thus the functional is defined up to the periods,
so we get a well-defined image in $\Jac(C)$.

\subsubsection*{\bf C}
The resulting map $\alpha:\Pic_0(C)\to\Jac(C)$
is called {\em the Abel-Jacobi map}. The Abel-Jacobi Theorem
states that it is an isomorphism. This theorem holds in tropical
geometry. It can be used to induce the tropical structure to $\Pic_0(C)$
and thus to all other $\Pic_k(C)$. The map $C\times\dots\times C\to\Pic_k(C)$
is a tropical morphism.

\subsection{The Riemann-Roch Theorem}
\subsubsection*{\bf A}
A divisor $D\in\Div(C)$ is called {\em effective} (we write $D\ge 0$)
if it is presented as $\sum a_jx_j$ with $a_j\ge 0$. If two effective
divisors $D,D'$ are linearly equivalent then they can be connected
by a deformation. Indeed we can find a rational function $\phi$
such that $D$ is the divisor of its zeroes and $D'$ is the divisor
of its poles. The full graph of $\phi$ is a tropical curve in
$C\times\T\PP^1$ (considered with $(x,y)$-coordinates).
Functions $y+t$, $t\in\T^\times$, define the required family of divisors
connecting $D$ and $D'$.

\subsubsection*{\bf B}
All effective divisor in a given equivalence class thus form
a connected space. It is denoted with $|D|$. In complex geometry
this space is the complex projective space $\C\PP^n$, where $n$
is $\dim H_0(D,\OO)$. In tropical geometry the space $|D|$ is also
kind of projective, though usually it is not isomorphic to $\T\PP^n$
for any $n$.

\subsubsection*{\bf C}
There is a tropical vector space, i.e. a space that contains $\T$ and
admits addition and multiplication by scalars (elements of $\T$) such that
$|D|$ comes as its tropical projectivisation. Namely, choose $D_0\in |D|$
and consider the space $V$
of all rational functions $\phi$ such that at every $x\in D$
either $\phi$ is regular or has a pole of order not more than the multiplicity
of $D_0$ at $x$. Clearly, if $\phi,\psi\in V$ then
$``a\phi+b\psi"\in V$ for any $a,b\in\T^\times$. Furthermore, each $D\in |D|$
corresponds to a ray in $V$, i.e. an element of $V$ up to a multiplication by
a scalar.

\subsubsection*{\bf D}
In a contrast with the complex case there are many tropical vector
space not isomorphic to $\T^n$. Nevertheless, we can define the dimension
of $|D|$ in the following way: we say that $\dim|D|\ge k$ if for any effective
divisor $D'$ of degree $k$ the space $|D-D'|$ is nonempty. The Riemann-Roch
theorem states that $\dim|D|-\dim|K-D|=d-g+1$.
Here $d=\deg D$, $g$ is the genus of the curve and $K$ is the divisor obtained
by taking every vertex of $C$ with the coefficient equal to its valence minus two.
The equivalene class of the divisor $K$ is called {\em the canonical class}.

\section{Tropical Gromov-Witten theory and enumerative geometry}

\subsection{Moduli spaces of tropical curves}
\subsubsection*{\bf A}
Consider all tropical curves of a given genus $g$ with $k$
distinct marked points up to equivalence generated by the tropical modifications.
If $g=0$ we assume in addition that $k\ge 3$ and if $g=1$ that $k>0$. As above
applying the tropical modifications and the operations inverse to them
we may assume that all marked points are at the ends of the curve.

\subsubsection*{\bf B}
If $g=0$ all such curves form a tropical variety of dimension $k-3$
denoted with $\MM_{0,k}$.
This tropical variety can be embedded to $\R^N$ for some $N$ as a balanced
complex by the so-called cross-ratio functions. A cross-ratio on $C$ is defined
by a choice of two ordered disjoint pairs of elements of the marking set $\{1,\dots,k\}$.
Namely, we connect the points in each pair with an embedded path (unique
in the tree $C$) and measure the length $l$ of the intersection. The corresponding
cross-ratio is $l$ if the orientations of the paths agree and $-l$ if they disagree.

\subsubsection*{\bf C}
If $g>0$ the space $\MM_{g,k}$ is not a tropical variety. The reason
is possible symmetries of the corresponding curves. To exclude
such symmetries one may introduce the {\em homotopy marking},
or a homotopy equivalence $h: F_g\to C$, where $F_g$ is the (topological)
wedge of $g$ copies of the circle $S^1$. The tropical curves enhanced
with such marking form the universal covering $\tilde\MM_{g,k}$
of $\MM_{g,k}$. The deck group is isomorphic to the group of outer automorphisms
of the free group on $g$ generators.
In Geometric Group Theory the space $\tilde\MM_{g,0}$
is known as the Outer Space.

\subsubsection*{\bf D}
The space $\tilde\MM_{g,k}$ has a natural local structure
of a tropical variety of dimension $3g-3+k$.
To verify this we consider a curve $C\in\tilde\MM_{g,k}$.
Near any vertex of valence $v \ge 4$ the curve $C$ maybe perturbed
and the perturbation space is parameterized by the tropical
variety $\MM_{0,v}$. Furthermore for each inner edge we can locally
deform its length. Thus near $C$ the space $\tilde\MM_{g,k}$
coincides with a product of tropical varieties. The dimension
$3g-3+k$ is the number of inner edges in the case when all inner
vertices of $C$ are 3-valent.

\subsubsection*{\bf E}
Note that $\tilde\MM_{g,k}$ is not compact and even does not
have a finite type. Nevertheless its quotient by the group of deck
transformation $\MM_{g,k}$ does have a finite type (in some sense)
and can be viewed as a {\em tropical orbifold}.

\subsection{Stable tropical curves}
\subsubsection*{\bf A}
As in the classical case we may compactify $\MM_{g,k}$.
To do this it is sufficient to allow the inner edges to assume
infinite length. The result is denoted with $\overline\MM_{g,k}$.
It is clearly compact. Furthermore, it is a tropical variety
as the locally the compactification construction just corresponds
to passing from $\R^N$ to $\T^N$ by taking the closure.

\subsubsection*{\bf B}
The elements of $\overline\MM_{g,k}$ are called {\em stable curves}.
To establish the correspondence with
the classical construction of the Deligne-Mumford
compactification by curves with several components it suffices to contract
all edges of infinite length (note that then we may have several intersecting
and self-intersecting components).

\subsubsection*{\bf C}
Let $X$ be a tropical variety. We may consider the space
$\MM_{g,k}(X)$ of all tropical morphisms from $C$ to $X$,
where $C\in\overline\MM_{g,k}$.
The image of $C$ is a 1-cycle in $X$. We may specify the homology
class of such 1-cycle.
To simplify our considerations
we assume for the rest of these lectures that $X=\T\PP^n$.
Then specifying the homology class is the same as fixing
the degree of the curve. We denote the corresponding
set of tropical morphisms with $\MM_{g,k,d}(\T\PP^n)$.

\subsubsection*{\bf D}
Clearly, $\MM_{g,k,d}(X)$ is a polyhedral complex (though not necessarily
of a pure dimension and its absence does not allow us to check the balancing
condition in general).
It follows from the tropical Riemann-Roch formula that
$$\dim\MM_{g,k,d}(\T\PP^n)\ge (n+1)d+(n-3)(1-g).$$

\subsubsection*{\bf E}
Let us consider a the closure in $\MM_{g,k,d}(\T\PP^n)$
tropical morphisms that are realizable by topological immersions
that locally deform in exactly $((n+1)d+(n-3)(1-g))$-dimensional space.
This subspace is called the {\em moduli space of regular curves},
we denote it with $\MM^{reg}_{g,k,d}(\T\PP^n)$.
It can be checked that it is a tropical variety
(of dimension $(n+1)d+(n-3)(1-g)$).
The elements of $\MM_{g,k,d}(\T\PP^n)\setminus \MM^{reg}_{g,k,d}(\T\PP^n)$
are called the {\em superabundant} curves.

\subsubsection*{\bf F}
As in the complex case it is possible to compactify $\MM^{reg}_{g,k,d}(\T\PP^n)$
to a compact tropical variety $\overline\MM^{reg}_{g,k,d}(\T\PP^n)$.
To do this we consider the morphisms from stable curves with the Kontsevich
condition on each component (of the complement of infinite length edges)
that is contracted to zero.

\subsubsection*{\bf G}
There are two special cases when there are no superabundant curves.
Namely, for immersions in the case when $n=2$ and in the case when $g=0$.

\subsection{Tropical Gromov-Witten theory}
\subsubsection*{\bf A}
Suppose that $\overline\MM_{g,k,d}(X)$ is a compact tropical variety
of dimension prescribed by the Riemann-Roch formula. Then
we have the evaluation map
$$\ev_j:\overline\MM_{g,k,d}(X)\to X$$ for each of the $k$
marked points. This map is a tropical morphism.

\subsubsection*{\bf B}
Given a tropical cycle $A\subset X$ we may pull it back by $\ev_j$ to
a tropical cycle $\ev^*_j(A)\subset\overline\MM_{g,k,d}(X)$. This
pull-back is especially easy to define in the case when $A$ is a Cartier
divisor or a transverse intersection of such divisors as we may just
pull-back the corresponding regular functions.

\subsubsection*{\bf C}
This gives an alternative (with respect to the usual tropical
stable intersection in $X$) product operation for the cycles in $X$.
First we may pull-back the cycles to $\overline\MM_{g,k,d}(X)$
and then take their tropical (stable) intersection in $\overline\MM_{g,k,d}(X)$.
These products are known as the Gromov-Witten numbers.

\subsubsection*{\bf D}
In the case of $g=0$ the resulting product in the clsssical (complex) set-up
is known as the {\em quantum} product. Its associativity is related to the so-called
WDVV equation. This equation holds in the tropical case as well.
As in the classical case it follows from the linear equivalence of all points
in $\overline\MM_{0,4}$ and, in particular, its boundary (sedentarity 1) points.

\subsubsection*{\bf E}
Geometrically, the Gromov-Witten numbers may be inerpreted as the
number of curves of genus $g$ and degree $d$ passing via the corresponding
$k$ cycles. Clearly, to have the non-zero number a certain condition on the
sum of the codimensions of the cycles has to hold.

\subsubsection*{\bf F}
It can be shown that in the regular case (in particular,
in the case $n=2$ or $g=0$ there is a correspondence between tropical
curves and their classical (real and complex) counterparts.
This gives a way to reduce computation of classical invariants
to finding the tropical curves and computing their real and complex multiplicities
(which are different, in general).
There are different algorithms for doing that. In the lectures we consider many
examples of tropical computations.

\subsubsection*{\bf G}
One of the powerful technique to enumerate tropical curves in $\T\PP^n$
is via the so-called {\em floor decompositions}. The technique is based
on the observation that each component of the complement of
vertical edges for a tropical curve in $\R^n$ is projected (vertically)
to a tropical curve in $\R^{n-1}$. This allows to find curves inductively
by dimension. Each such component is called a {\em floor}.
The floor diagram is the graph obtained by contracting each floor to a point.
The resulting computation is especially simple in the case of $n=2$.

\section{Historical remarks and credits}
Passing to the equivalent of tropical limit for counting holomorphic curves
and holomorphic disks respectively was suggested by Kontsevich and Fukaya
around 2000. Formalization of the limiting objects as tropical varieties and
their geometry as Tropical Geometry were introduced by Mikhalkin and
Sturmfels in 2002. Applications of Tropical Geometry in classical
enumerative (complex and real) curves were started by Mikhalkin in 2002.
These applications can be considered as a build-up on patchworking
(Viro, 1979) and amoebas (Gelfand, Kapranov, Zelevinski, 1994),
in particular, non-Archimedean amoebas (Kapranov, 2000).
A general understanding of passing from classical to tropical Mathematics
as a sort of {\em dequantization} goes back to Maslov and his school
(Kolokoltsov, Litvinov et al.) in the early 90s.
Below are some references for the material of these lectures.


\begin{thebibliography}{99}
\bibitem{} E. Brugall\'e, G. Mikhalkin, {\em Enumeration of curves via floor diagrams},
http://arxiv.org/abs/0706.0083.
\bibitem{} A. Gathmann, H. Markwig, {\em Kontsevich's formula and
the WDVV equations in tropical geometry}, http://arxiv.org/abs/math.AG/0509628.
\bibitem{} G. Mikhalkin, {\em Amoebas of algebraic varieties and tropical geometry,}
http://arxiv.org/abs/math.AG/0403015, in {\em Different Faces of Geometry},
International Mathematical Series , Vol.  3
(Donaldson, Simon; Eliashberg, Yakov; Gromov, Mikhael (Eds.))
2004, 257--300.
\bibitem{} G. Mikhalkin, {\em Tropical geometry and its application},
Proceedings of the ICM 2006 Madrid, 25 pages,
http://arxiv.org/abs/math.AG/0601041.
\bibitem{} G. Mikhalkin, {\em Tropical Geometry},
ongoing book project, partial material available at http://www.math.toronto.edu/mikha/book.ps.
\bibitem{} G. Mikhalkin, I. Zharkov,
{\em Tropical curves, their Jacobians and theta-functions},
http://arxiv.org/abs/math.AG/0612267.
\bibitem{} J. Richter-Gebert, B. Sturmfels, Th. Theobald,
{\em First steps in tropical geometry}, http://arxiv.org/abs/math.AG/0306366.
\bibitem{} M. Vidgeland, {\em Smooth tropical surfaces with infinitely many tropical lines},
http://arxiv.org/abs/math/0703682.
\bibitem{} O. Viro, {\em Dequantization of real algebraic geometry on logarithmic paper},
Proceedings of the 3d ECM, Barcelona 2000. http://arxiv.org/abs/math.AG/0005163.
\end{thebibliography}
\end{document}